\documentclass{article}
\usepackage{amssymb}
\usepackage{latexsym}
\newcommand{\R}{{\mathbb R}}
\newcommand{\N}{{\mathbb N}}

\newcommand{\bx}{\hfill{$\Box $ }}

\newtheorem{lemma}{Lemma}[section]
\newtheorem{corollary}{Corollary}[section]
\newtheorem{proposition}{Proposition}[section]
\newtheorem{theorem}{Theorem}[section]
\title{On Krawtchouk polynomials}
\author{R.Coleman\\Laboratoire LJK-Universit\'e de Grenoble,\\ Tour-IRMA,\\51, rue des Math\'ematiques,\\Domaine universitaire de Saint-Martin-d'H\`eres, France.}
\begin{document}
\maketitle

Krawtchouk polynomials play an important role in coding theory and are also
useful in graph theory and number theory. Although the basic properties of 
these polynomials are known to some extent, there is, to my knowledge, no 
detailed development available. My aim in writing this article is to fill in this gap.\\

\noindent {\bf Notation} In the following we will use capital letters for
(algebraic) polynomials, for example $P$ or $P(X)$; for the 
polynomial function associated with a polynomial $P$, we will use small 
letters in the parentheses, for example $P(x)$.\\

\section{Definition and first properties}
To begin with, we define a collection of real polynomials $P_j$, for $j\geq
0$, which we will use in the definition of Krawtchouk polynomials. We set 
$P_0=1$ and, for $j\geq 1$, 
$$
P_j(X) = \frac{1}{j!}X(X-1)\ldots (X-j+1).
$$
Clearly, for $i\geq j\geq 0$, $P_j(i)={i\choose j}$. Now we take $n,m,s\in \N$, such
that $m\leq n$ and $s\geq 2$, and set
$$
K_{m,n,s}(X)=\sum _{j=0}^m(-1)^jP_j(X)P_{m-j}(n-X)(s-1)^{m-j}.
$$
The real polynomial $K_{m,n,s}$ is called a Krawtchouk polynomial.  For 
every fixed pair $(n,s)\in \N\times \N\setminus \{0,1\}$, there
is a set of $n+1$ Krawtchouk polynomials $K_{m,n,s}$. Clearly $K_{0,n,s}$ is the
constant polynomial $1$. To simplify the notation, when the indices $n$ and $s$ 
are understood, we will write $K_m$ for $K_{m,n,s}$.

\begin{proposition} $K_m$ is a polynomial of degree $m$, whose leading
coefficient is $\frac{-s^m}{m!}$.
\end{proposition}

\noindent\textsc{proof} We first notice that
$$
P_j(X)P_{m-j}(n-X) = \frac{(-1)^{m-j}}{j!(m-j)!}X^m
\;+\;\textrm{terms of degree smaller than $m$},
$$
which implies that 
$$
K_m(X) = \sum _{j=0}^m\frac{(-1)^m}{j!(m-j)!}(s-1)^{m-j}X^m 
\;+\;\textrm{terms of degree smaller than $m$}.
$$
However
$$
\sum _{j=0}^m\frac{(-1)^m}{j!(m-j)!}(s-1)^{m-j} =
\frac{(-1)^m}{m!}\sum _{j=0}^m{m\choose j}(s-1)^{m-j} = \frac{(-s)^m}{m!},
$$ 
and hence the result. \bx\\

We have a simple expression for the value of the polynomial $K_{m,n,s}$ at certain
integers.

\begin{proposition}
For $i\in\{0,\ldots ,n\}$, $K_{m,n,s}(i)$ is the coefficient of $X^m$ in the
product
$$
(1+(s-1)X)^{n-i}(1-X)^i.
$$
\end{proposition}

\noindent\textsc{proof} First we have
$$
K_{m,n,s}(0) =\sum _{j=0}^m(-1)^jP_j(0)P_{m-j}(n-0)(s-1)^{m-j} =
P_m(n)(s-1)^m = {n\choose m}(s-1)^m,
$$
which is the coefficient of $X^m$ in the expression $(1+(s-1)X)^n$. It is also 
easy to see that
$$
K_{m,n,s}(n) = (-1)^mP_m(n) = (-1)^m{n\choose m},
$$
which is the coefficient of $X^m$ in the expression $(1-X)^n$. So the result is
true for $i=0$ and $i=n$.

Now let us consider the case where $1\leq i\leq n-1$. We have
$$
(1+(s-1)X)^{n-i} = 1+{n-i\choose 1}(s-1)X+{n-i\choose 2}(s-1)^2X^2+\ldots +
(s-1)^{n-i}X^{n-i}
$$
and 
$$
(1-X)^i = 1-{i\choose 1}X+{i\choose 2}X^2 -\ldots + (-1)^iX^i.
$$
It follows that the coefficient of $X^m$ in the product of $(1+(s+1)X)^{n-1}$ 
and $(1-X)^i$ is 
$$
\sum _{j=0}^m(-1)^j{i\choose j}{n-i\choose m-j}(s-1)^{m-j}=K_{m,n,s}(i).
$$
This ends the proof.\bx\\

\noindent {\bf Remark} Further on we will use the expressions for $K_{m,n,s}(0)$
and $K_{m,n,s}(n)$ found here, namely
$$
K_{m,n,s}(0)={n\choose m}(s-1)^m \quad \textrm{and}\quad 
K_{m,n,s}(n)=(-1)^m{n\choose m}.
$$

\begin{corollary} \label{cor1a} The following statement follows directly from
the proposition:
$$
K_{1,n,s}(i)= (n-i)(s-1)-i=(s-i)n-si.
$$
\end{corollary}

It is very easy to find an explicit expression for the polynomials $K_{1,n,s}$.

\begin{proposition} \label{prop1aa} We have
$$
K_{1,n,s}(X) = n(s-1)-sX.
$$
\end{proposition}

\noindent \textsc{proof} Using Corollary \ref{cor1a}, we obtain 
$K_{1,n,s}(0)=(s-1)n$ and $K_{1,n,s}(1)=(s-1)n-s$, from which we deduce the 
two coefficients of $K_{1,n,s}$.\bx\\

From what we have seen, we may obtain a useful recurrence relation.

\begin{proposition} \label{prop1a} For $1\leq m,i\leq n$, we have
$$
K_{m,n,s}(i) = K_{m,n,s}(i-1)-K_{m-1,n,s}(i-1)-(s-1)K_{m-1,n,s}(i).
$$
\end{proposition}

\noindent \textsc{proof} We first notice the identity
$$
(1+(s-1)X)^{n-i}(1-X)^i(1+(s-1)X)=
\left((1+(s-1)X)^{n-(i-1)}(1-X)^{i-1}\right)(1-X).
$$
The coefficient of $X^m$ on the left-hand side is
$$
K_{m,n,s}(i)+(s-1)K_{m-1,n,s}(i)
$$
and on the right-hand side is
$$
K_{m,n,s}(i-1)-K_{m-1,n,s}(i-1).
$$
This gives us the result.\bx

\begin{corollary} \label{cor1b} We have
$$
K_{m,n,s}(X) = K_{m,n,s}(X-1)-K_{m-1,n,s}(X-1)-(s-1)K_{m-1,n,s}(X).
$$
\end{corollary}

\noindent \textsc{proof} The polynomials $K_{m,n,s}(X)$ and $K_{m,n,s}(X-1)$
have the same degree and leading coefficient, so the degree of their 
difference is at most $m-1$. As the polynomials
$$
K_{m,n,s}(X)-K_{m,n,s}(X-1)\quad \textrm{and}\quad
-(K_{m-1,n,s}(X-1)+(s-1)K_{m-1,n,s}(X)) 
$$
have the same value at $n$ points and $m-1\leq n$, these polynomials are the same and
we have the result we were looking for.\bx\\

\noindent {\bf Remark} Using the recurrence relation of the proposition, we may 
successively calculate
$K_{m,n,s}(1)$, $K_{m,n,s}(2)$ and so on up to $K_{m,n,s}(m)$ and so determine 
the polynomial $K_{m,n,s}$.\\
 
\noindent {\bf Example} To calculate $K_{2,n,s}$, it is sufficient to obtain
$K_{2,n,s}(0)$, $K_{2,n,s}(1)$ and $K_{2,n,s}(2)$. To find $K_{2,n,s}(0)$, we 
can use the the expression for $K_{m,n,s}(0)$ found above. We can calculate 
$K_{2,n,s}(1)$ from the recurrence relation of the proposition: we need 
$K_{2,n,s}(0)$ (already calculated), $K_{1,n,s}(0)$ and $K_{1,n,s}(1)$, which 
can be found from the expression for $K_{1,n,s}$ (Proposition \ref{prop1aa}). We 
may also determine $K_{2,n,s}(2)$ from the recurrence relation of the
proposition: we need $K_{2,n,s}(1)$ (already calculated) and $K_{1,n,s}(1)$ 
and $K_{1,n,s}(2)$, which can be found using the expression for $K_{1,n,s}$.
 
\section{A summation formula}
The recurrence relation found in the last section enables us to find a formula for
the sum of successive Krawtchouk polynomials. We need a preliminary result.

\begin{lemma}\label{lemma2a}
For $j\geq 0$, we have
$$
P_j(X) + P_{j+1}(X) = P_{j+1}(X+1).
$$
\end{lemma}

\noindent \textsc{proof} Let 
$$
Q(X) = P_j(X) + P_{j+1}(X) - P_{j+1}(X+1).
$$
For any $i\geq j+1$, we have
$$
Q(i) = P_j(i) + P_{j+1}(i) - P_{j+1}(i+1) = {i\choose j} + {i\choose j+1} -
{i+1\choose j+1} = 0.
$$
Therefore $Q=0$.\bx\\

We now turn to the summation formula.

\begin{theorem} For $1\leq m\leq n$, we have
$$
\sum _{k=0}^mK_{k,n,s}(X) = K_{m,n-1,s}(X-1).
$$
\end{theorem}

\noindent \textsc{proof} We will prove this result by induction on $m$. First
$$
K_{0,n,s}(X) + K_{1,n,s}(X) = 1 + n(s-1)-sX = (n-1)(s-1)-s(X-1) = K_{1,n-1,s}(X-1),
$$
therefore the theorem is correct for $m=1$. 

Suppose now that the result is true for
$m-1$ and consider the case $m$. Using Corollary \ref{cor1b}, we
obtain
$$
\sum _{k=1}^mK_{k,n,s}(X) =
K_{m,n,s}(X-1)-K_{0,n,s}(X-1)-(s-1)\sum_{k=0}^{m-1}K_{k,n,s}(X),
$$
which implies that 
$$
\sum _{k=0}^mK_{k,n,s}(X) = K_{m,n,s}(X-1)-(s-1)\sum_{k=0}^{m-1}K_{k,n,s}(X),
$$
because $K_{0,n,s}(X-1)=1=K_{0,n,s}(X)$. Now
$$
K_{m,n,s}(X-1) = \sum _{j=0}^m(-1)^jP_j(X-1)P_{m-j}(n-X+1)(s-1)^{m-j}
$$
and, by hypothesis, 
\begin{eqnarray*}
(s-1)\sum_{k=0}^{m-1}K_{k,n,s}(X) &=& (s-1)K_{m-1,n-1,s}(X-1)\\
		&=& (s-1)\sum _{j=0}^{m-1}(-1)^jP_j(X-1)P_{m-1-j}((n-1)-(X-1))(s-1)^{m-1-j}\\
		&=& \sum _{j=0}^{m-1}(-1)^jP_j(X-1)P_{m-1-j}(n-X)(s-1)^{m-j}.
\end{eqnarray*}
Using Lemma \ref{lemma2a}, we obtain
\begin{eqnarray*}
\sum _{k=0}^mK_{k,n,s}(X) &=& \sum _{j=0}^{m-1}(-1)^jP_j(X-1)P_{m-j}(n-X)(s-1)^{m-j} +
(-1)^mP_m(X-1)\\
	&=& \sum _{j=0}^m(-1)^jP_j(X-1)P_{m-j}(n-X)(s-1)^{m-j}\\
	&=& K_{m,n-1,s}(X-1).
\end{eqnarray*}
Therefore the result is true for $m$. This finishes the induction step.\bx\\
	
\section{Inner products}
For $n\in \N$, let us write $\R _n[X]$ for the set of real polynomials of 
degree not greater than $n$. $\R _n[X]$ is a real vector space of 
dimension $n+1$. Fixing $s\geq 2$, we define an inner product $\langle \cdot, \cdot
\rangle$ on $\R _n[X]$ in the following way:
$$
\langle A,B\rangle = \sum _{i=0}^n{n\choose i}(s-1)^iA(i)B(i).
$$
If $0\leq m\leq n$, then the polynomial $K_{m,n,s}\in \R _n[X]$. In the next
proposition, we will drop the second and third parameters to simplify the
notation.

\begin{proposition} For $0\leq k,l\leq n$, we have
$$
\langle K_k,K_l\rangle = 
\cases{0&$k\neq l$\cr
s^n(s-1)^k{n\choose k}&$k=l$\cr}
$$
\end{proposition}

\noindent \textsc{proof} Let us consider the real polynomial in two
variables
$$
P(X,Y) = \sum _{k=0}^n\sum _{l=0}^n\left(\sum _{i=0}^n{n\choose
i}(s-1)^iK_k(i)K_l(i)\right)X^kY^l. 
$$
Then
\begin{eqnarray*}
P(X,Y) &=& \sum _{i=0}^n{n\choose i}(s-1)^i\sum _{k=0}^nK_k(i)X^k
\sum _{l=0}^nK_l(i)Y^l\\
       &=& \sum _{i=0}^n{n\choose i}(s-1)^i(1+(s-1)X)^{n-i}(1-X)^i
       (1+(s-1)Y)^{n-i}(1-Y)^i\\
       &=& \sum _{i=0}^n{n\choose
       i}\left((s-1)(1-X)(1-Y)\right)^i\left((1+(s-1)X)(1+(s-1)Y)\right)^{n-i}\\
       &=& \left((s-1)(1-X)(1-Y)+(1+(s-1)X)(1+(s-1)Y)\right)^n.
\end{eqnarray*}
After simplification, we obtain
$$
P(X,Y) = s^n(1+(s-1)XY)^n = \sum _{i=0}^ns^n{n\choose i}(s-1)^iX^iY^i.
$$
Therefore, if $k\neq l$, then
$$
\sum _{i=0}^n{n\choose i}(s-1)^iK_k(i)K_l(i) = 0
$$
and, if $k=l$, then
$$
\sum _{i=0}^n{n\choose i}(s-1)^iK_k(i)K_l(i) = s^n{n\choose k}(s-1)^k.
$$
This ends the proof.\bx\\

\begin{corollary} \label{cor2a} If $A$ is a real polynomial such that $\deg (A)=d<m$, then
$$
\langle K_m,A\rangle = 0.
$$
\end{corollary} 

\noindent \textsc{proof} The polynomials $K_k$, $0\leq k\leq n$, form a basis
of $\R _n[X]$. If $\deg (A)=d$, then $A=\sum _{k=0}^d\lambda _kK_k$ and
$$
\langle K_m,A\rangle = \sum _{k=0}^d\lambda _k\langle K_m , K_k\rangle =0.
$$
This ends the proof.\bx

\section{Another recurrence relation} 
We have already seen one recurrence relation involving Krawtchouk polynomials. 
In this section we present another such relation, which will be useful when
studying the roots of Krawtchouk polynomials. Once again, we will drop the
parameters $n$ and $s$ from $K_{m,n,s}$ to simplify the notation.

\begin{theorem} \label{th3a} For $1\leq m\leq n-1$, we have the relation
$$
(m+1)K_{m+1} = (m+(s-1)(n-m)-sX)K_m-(s-1)(n-m+1)K_{m-1}.
$$
\end{theorem} 

\noindent \textsc{proof} We have seen that
$$
\sum _{m=0}^nK_m(i)X^m = (1+(s-1)X)^{n-i}(1-X)^i.
$$
If we differentiate both sides of the equation, then we obtain
$$
\sum _{m=0}^{n-1}(m+1)K_{m+1}(i)X^m =
(s-1)(n-i)(1+(s-1)X)^{n-i-1}(1-X)^i-i(1+(s-1)X)^{n-i}(1-X)^{i-1}.
$$
We now multiply both sides of this expression by $(1+(s-1)X)(1-X)$. Writing $A$
and $B$ for the left- and right-hand sides, we have
$$
(1+(s-1)X)(1-X)A = \sum _{m=0}^{n-1}(m+1)K_{m+1}(i)X^m +(s-2)\sum
_{m=1}^nmK_m(i)X^m - (s-1)\sum _{m=2}^{n+1}(m-1)K_{m-1}(i)X^m
$$
and     
\begin{eqnarray*}
(1+(s-1)X)(1-X)B &=& \left((1-X)(s-1)(n-i)-i(1+(s-1)X)\right)
		(1+(s-1)X)^{n-i}(1-X)^i\\
		 &=& ((sn-n-is)+(n-sn)X)\sum _{m=0}^nK_m(i)X^m\\
		 &=& (sn-n-is)\sum _{m=0}^nK_m(i)X^m+(n-sn)\sum
		 _{m=1}^{n+1}K_{m-1}(i)X^m.
\end{eqnarray*}
Taking the difference of $(1+(s-1)X)(1-X)A$ and $(1+(s-1)X)(1-X)B$, we obtain
$$
\sum _{m=0}^{n-1}(m+1)K_{m+1}(i)X^m- \sum _{m=0}^n(m+(s-1)(n-m)-si)K_m(i)X^m +
\sum _{m=1}^n(s-1)(n-m+1)K_{m-1}(i)X^m = 0.
$$
Hence, for $1\leq m\leq n-1$ and $0\leq i\leq n$,
$$
(m+1)K_{m+1}(i) - (m+(s-1)(n-m)-si)K_m(i) + (s-1)(n-m+1)K_{m-1}(i)=0.
$$
However, the polynomial
$$
C_m = (m+1)K_{m+1} - (m+(s-1)(n-m)-sX)K_m + (s-1)(n-m+1)K_{m-1}
$$
is of degree not greater than $n$ and has $n+1$ roots. It follows that $C_m=0$,
which ends the proof.\bx\\ 

\noindent {\bf Remark} This recurrence relation also allows us to find the
polynomials $K_m$ successively. For example, $K_2$ can be obtained from $K_1$
and $K_0$, which we already know. After some calculation, we find
$$
K_2 = \frac{1}{2}\left((s-1)^2n(n-1)-s(2ns-2n-s+2)X+s^2X^2\right).
$$
When $s=2$, the expression is much simpler:
$$
K_2 = \frac{(n-2X)^2-n}{2}.
$$

\section{Roots of Krawtchouk polynomials}
In this section we will see that all the roots of a Krawtchouk polynomial are
real and distinct. We obtain an interesting relation between the roots of
successive Krawtchouk polynomials. To simplify the notation, we set
$$
a_m = m+(s-1)(n-m) \qquad \textrm{and} \qquad b_m = (s-1)(n-m+1).
$$
Notice that $a_m>0$ and $b_m>0$.\\

\begin{proposition} For $1\leq m\leq n$, the polynomial $K_m$ has $m$ 
distinct real roots in the interval $(0,n)$.
\end{proposition}

\noindent \textsc{proof} As
$$
0 = \langle K_0,K_m\rangle = \sum _{i=0}^n{n\choose i}(s-1)^iK_m(i)
$$
and $K_m(0)>0$, the real-valued polynomial function $K_m(x)$ changes sign in the 
interval $(0,n)$. Suppose that $K_m(x)$ changes sign at the points $x_1<\ldots
<x_d$. These points are roots of $K_m$. As $\deg (K_m)=m$, we have $d\leq m$. 
Suppose that $d<m$ and let
$$
S= \prod _{i=1}^d(X-x_i).
$$
The polynomial $S$ has $x_1,\ldots ,x_d$ as roots. As these roots are simple,
the derivative $S'$ does not have any of these points as roots. It follows
that the real-valued polynomial function $S(x)$ changes sign at 
$x_1,\ldots ,x_d$. Thus the polynomial functions $K_m(x)$ and $S(x)$ change
signs at the same points. Hence the product $S(x)K_m(x)$ is strictly 
positive (or strictly negative), except at the points $x_1, \ldots ,x_d$, 
where its value is $0$. Therefore $\langle S,K_m\rangle > 0$ 
(or $<0$). However, this contradicts Corollary \ref{cor2a} and so $d=m$, i.e. 
$K_m$ has $m$ distinct roots in the interval $(0,n)$.\bx\\

\noindent {\bf Remark} As $\deg (K_m)=m$, the roots of $K_m$ are necessarily
simple.\\

We will now look at the `interlacing' property of the roots of Krawtchouk
polynomials.\\

\noindent \begin{theorem} If the roots of $K_m$ are $x_1<\ldots <x_m$ and
those of $K_{m+1}$ are $y_1<\ldots <y_{m+1}$, then
$$
0<y_1<x_1<y_2<x_2<\ldots <x_m<y_{m+1}<n.
$$
\end{theorem}

\noindent \textsc{proof} We will prove the result by induction on $m$. We first
notice that, from Theorem \ref{th3a}, we have
$$
(m+1)K_{m+1} = (a_m-sX)K_m-b_mK_{m-1}.
$$
We begin with the case $m=1$. Let $x_1$ be the unique root of $K_1$. As 
$K_0=1$ and $b_1>0$, $K_2(x_1)<0$. Also, $K_2(0)=(s-1)^2{n\choose 2}>0$ and 
$K_2(n)=(-1)^2{n\choose 2}>0$, and so
$K_2$ has a root in each of the intervals $(0,x_1)$ and $(x_1,n)$. This
proves the result for $m=1$.

Suppose now that the property is true for $m-1$ and consider the case 
$m$. Let $x_1<\ldots <x_m$ be the roots of $K_m$. Then $K_{m-1}(x_1)>0$,
$K_{m-1}(x_2)<0$, $\ldots$ This implies that $K_{m+1}(x_1)<0$,
$K_{m+1}(x_2)>0$, $\ldots$ Therefore $K_{m+1}$ has a root in the interval
$(x_i,x_{i+1})$, for $i\in\{1,\ldots ,m-1\}$. Also, $K_{m+1}(0)>0$ and
$K_{m+1}(x_1)<0$, therefore $K_{m+1}$ has a root in the in the interval 
$(0,x_1)$. We claim that $K_{m+1}$ has a root in the interval $(x_m,n)$. We
observe that $K_{m+1}(n)=(-1)^{m+1}{n\choose m+1}$ and that $K_{m+1}(x_m)>0$, if
$m$ is even, and $K_{m+1}(x_m)<0$, if $m$ is odd. This implies that $K_{m+1}$
has a root in the interval $(x_m,n)$. Thus the result is true for $m$. This
finishes the proof.\bx\\ 

In the case where $s=2$, we can say a little more about the roots. This follows
from a simple symmetry relation.

\begin{lemma} For $x\in \R$, we have
$$
K_{m,n,2}(n-x) = (-1)^mK_{m,n,2}(x).
$$
\end{lemma}

\noindent \textsc{proof} We have
\begin{eqnarray*}
K_{m,n,2}(n-x)&=&\sum _{j=0}^m(-1)^jP_j(n-x)P_{m-j}(x)\\
	&=&\sum _{t=0}^m(-1)^{m-t}P_{m-t}(n-x)P_t(x)\\
	&=&(-1)^m\sum _{t=0}^m(-1)^tP_{m-t}(n-x)P_t(x)\\ 	
 	&=&(-1)^mK_{m,n,2}(x).
\end{eqnarray*}
This ends the proof.\bx\\

\begin{proposition} The roots of $K_{m,n,2}$ are symmetric with respect to
$\frac{n}{2}$. In particular, if $m$ is odd, then $\frac{n}{2}$ is a root of 
$K_{m,n,2}$.
\end{proposition}

\noindent \textsc{proof} It is sufficient to notice that, if $x_1$ is a root,
then so is $n-x_1$.\bx\\

\noindent {\bf Remark} It is particularly interesting to notice that, if $m$ is
odd and $n$ even, then $K_{m,n,2}$ has an integer root, namely $\frac{n}{2}$.

\end{document}